\documentclass[11pt,letterpaper]{amsart}
\usepackage[T1]{fontenc}
\usepackage[utf8]{inputenc}
\usepackage{amsmath,amssymb,amsthm}
\usepackage[colorlinks=true,linkcolor=blue,citecolor=blue,urlcolor=blue]{hyperref}
\newcommand{\C}{\mathbb C}
\newcommand{\R}{\mathbb R}
\newcommand{\Q}{\mathbb Q}
\newcommand{\OO}{\mathcal O}
\newcommand{\m}{\mathfrak m}
\newtheorem{theorem}{Theorem}
\newtheorem{lemma}[theorem]{Lemma}
\newtheorem{corollary}[theorem]{Corollary}
\theoremstyle{remark}

\title[A geometric characterization of weighted homogeneity]{A geometric characterization of weighted homogeneity for reduced complex analytic germs}
\author{Diogo da Silva Machado}
\author{Jos\'e Seade}
\date{}

\begin{document}
\begin{abstract}
We prove that a reduced complex analytic germ is weighted homogeneous if and only if, in suitable local coordinates, it admits a holomorphic vector field vanishing at the origin and pointing strictly outward along a single sufficiently small link. On a connected link, transversality suffices. The criterion applies in arbitrary codimension, including to non-isolated singularities. The proof uses a compact composition operator to control the spectrum on the Zariski tangent space and obtains a positive grading from the resulting spectral condition.
\end{abstract}
\maketitle

\section{The characterization}

In weighted homogeneous coordinates, the Euler vector field points strictly outward along every sufficiently small round link. We prove a converse for reduced complex analytic germs. If, in some local coordinates, a holomorphic vector field vanishes at the origin and points strictly outward along a single sufficiently small link, then the germ is weighted homogeneous. On a connected link, transversality suffices.

A germ $(X,0)$ is \emph{weighted homogeneous} if it admits an embedding in $(\C^N,0)$ whose ideal is generated by polynomials homogeneous for positive integer weights.

An embedding is \emph{minimal} if its ideal satisfies $I\subset\m^2$, where $\m=(z_1,\ldots,z_N)\subset\OO_{\C^N,0}$. A holomorphic vector field on $(X,0)$ is a derivation $v$ of $\OO_{X,0}$. It vanishes at $0$ if $v(\m_X)\subset\m_X$, where $\m_X$ is the maximal ideal of $\OO_{X,0}$. In that case its linear part is the induced operator on $\m_X/\m_X^2$. In a minimal embedding its spectrum agrees with that of the linear part of any ambient extension.

Put $\rho(z)=\sum_i|z_i|^2$ and $L_\varepsilon=X\cap\{\rho=\varepsilon^2\}$. All links below are taken in a sufficiently small representative of the germ. By \emph{transversal} we mean that
$d\rho(\operatorname{Re}v)\ne0$ on $L_\varepsilon$.
The field points strictly outward when this expression is positive. These conditions are imposed at every point of the link, including the points of $L_\varepsilon\cap\operatorname{Sing}X$. A field used on a link is understood to be defined near $X\cap\overline B_\varepsilon$.

\begin{theorem}\label{thm:main}
Let $(X,0)$ be a reduced complex analytic germ. The following are equivalent:
\begin{enumerate}
\item\label{item:wh} $(X,0)$ is weighted homogeneous.
\item\label{item:single} In suitable local coordinates, $X$ admits a holomorphic vector field vanishing at $0$ and pointing strictly outward along a single sufficiently small link.
\item\label{item:spectrum} $X$ admits a holomorphic vector field vanishing at $0$ whose linear spectrum lies in an open half-plane bounded by a line through $0$.
\end{enumerate}
If the link is connected, strict outward pointing in \eqref{item:single} can be replaced by transversality.
\end{theorem}

The connected-link formulation applies to every positive-dimensional irreducible germ. For example, the germ at the zero matrix of the variety of complex $3\times4$ matrices of rank at most two is irreducible of codimension two. Its ideal is minimally generated by the four maximal minors, and its singular locus consists of the matrices of rank at most one \cite[Chapter~1]{BV}. The ordinary Euler field is tangent to this cone and points strictly outward along its links. This gives a non-isolated example outside the class of complete intersections. The choice of coordinates matters even for isolated hypersurface singularities \cite[Proposition~8.2]{LZ2}.

This note springs from our preprint \cite{MS} on generic vector fields on isolated hypersurface germs. There, for isolated hypersurface singularities, we conjectured that weighted homogeneity is equivalent both to the existence of a holomorphic vector field transversal to the link and to that of one with a non-degenerate ambient extension; we proved the second criterion for even dimension and the first for odd dimension. The remaining cases were conjectured to hold for isolated hypersurface germs in general, and this was settled by Liu and Zhang \cite{LZ1,LZ2}. Here we prove the transversality criterion in full, and in much greater generality: Theorem~\ref{thm:main} holds for every reduced complex analytic germ, of arbitrary codimension and possibly with non-isolated singularities. It also characterizes weighted homogeneity by a spectral condition on the linear part of the field, which implies non-degeneracy. The argument is short and self-contained, and uses no defining equations of a particular form.

The proof separates into two steps. The first converts inward transversality into a spectral condition by a compact composition operator. The second converts this condition into a positive grading, using the formal semisimple part of a derivation. This places the criterion alongside the grading results of Scheja and Wiebe \cite{SW73,SW80} and Morvan's theorem that an analytic germ admitting a contracting automorphism is weighted homogeneous \cite{Morvan}. Indeed, the spectral condition produces such an automorphism by integrating a scalar multiple of the field. We give a direct proof to make the infinitesimal mechanism explicit.

\section{From transversality to the spectrum}

\begin{lemma}\label{lem:spectrum}
Let $(X,0)\subset(\C^N,0)$ be a reduced germ. If a holomorphic vector field $v$ vanishing at $0$ points strictly inward along a small link, every eigenvalue $\lambda$ of its linear part satisfies $\operatorname{Re}\lambda<0$.
\end{lemma}

\begin{proof}
Choose the representative so that every irreducible component of $W=X\cap B_\varepsilon$ contains $0$. The real-time flow $\varphi_t$ of $v$ is defined on $\overline W$ for $t\ge0$ and fixes $0$. Strict inward pointing at the boundary gives
\[
 \varphi(\overline W)\Subset W,\qquad \varphi=\varphi_1.
\]
We use the compact-open topology on $\OO(W)$, which is the canonical Fr\'echet topology because $W$ is reduced \cite[Chapter~V, \S6]{GR}. By Cartan's Theorem~B and the open mapping theorem, the restriction $\OO(B_\varepsilon)\to\OO(W)$ is a quotient map. Thus a sequence bounded on compact subsets of $W$ admits lifts uniformly bounded on any prescribed closed ball compactly contained in $B_\varepsilon$. Indeed, the image of the open unit ball for the corresponding supremum seminorm is a neighborhood of zero, and the sequence lies in a fixed multiple of that neighborhood.

Let $H=H^\infty(W)$, with its supremum norm, and let $\mathcal M\subset H$ consist of functions vanishing at $0$. To see that $T(u)=u\circ\varphi$ is compact, choose a closed ball inside $B_\varepsilon$ whose interior contains $\varphi(\overline W)$. For any bounded sequence in $H$, the preceding lifting argument and Montel's theorem on the interior of that ball give a subsequence converging uniformly on $\varphi(\overline W)$. The composed subsequence converges in $H$, which is complete in the supremum norm. The operator $T$ preserves $\mathcal M$.

Suppose $u\in\mathcal M$ and $Tu=cu$, with $|c|\ge1$. Then
\[
 \|u\|_W\le\max_{\varphi(\overline W)}|u|\le\|u\|_W.
\]
Thus $|u|$ attains its supremum in $W$. The maximum principle makes $u$ constant on an irreducible component containing a point of maximum. That component contains $0$, so the maximum is zero. Hence $u=0$, and $T|_{\mathcal M}$ has no eigenvalue of modulus at least one.

The first-jet map
\[
 j:\mathcal M\longrightarrow\m_X/\m_X^2
\]
is continuous because the map $\OO(W)\to\OO_{X,0}/\m_X^2$ descends from the continuous ambient jet map through the quotient map above. Restrictions of ambient linear forms supply every first jet, so $j$ is surjective. Its kernel $\mathcal N$ is closed and $T$-invariant. On $\mathcal M/\mathcal N$, the induced operator is the pullback of the linear flow and has eigenvalues $e^\lambda$.

Every nonzero eigenvalue $c$ of this quotient is an eigenvalue of $T|_{\mathcal M}$. Otherwise the Fredholm alternative would make $T-c$ surjective on $\mathcal M$, hence on the finite-dimensional quotient, contradicting the fact that $c$ is an eigenvalue there. Consequently $|e^\lambda|<1$, as required.
\end{proof}

\section{From the spectrum to a positive grading}

\begin{lemma}\label{lem:grading}
Let $(X,0)$ be a complex analytic germ, not necessarily reduced or isolated. If it admits a holomorphic vector field vanishing at $0$ whose linear spectrum lies in an open half-plane bounded by a line through $0$, then $(X,0)$ is weighted homogeneous.
\end{lemma}

\begin{proof}
Choose a minimal embedding with ideal $I\subset\OO_N=\C\{z_1,\ldots,z_N\}$ and an ambient extension $D$ of the field. Then $D(I)\subset I$. Work in the completion $\widehat\OO_N=\C[[z]]$, with ideal $\widehat I=I\widehat\OO_N$.

On each finite-dimensional algebra $\widehat\OO_N/\m^k$, the semisimple part of $D$ is again a derivation and is a polynomial in $D$, so it preserves the image of $\widehat I$. These semisimple parts commute with the quotient maps and define a formal derivation $S$ preserving $\widehat I$. Lifting an eigenbasis of $\m/\m^2$ through the successive quotients gives formal coordinates in which
\[
 S=\sum_{i=1}^N\lambda_i z_i\partial_{z_i},
\]
where the $\lambda_i$ are the eigenvalues of the original linear part. Here only the formal Chevalley decomposition is needed; see \cite{SW73,SW80}.

Choose $\theta\in\R$ and $\delta>0$ such that
\[
 \operatorname{Re}(e^{i\theta}\lambda_i)\ge\delta
 \qquad(1\le i\le N).
\]
For a monomial $z^a$ this gives
\[
 \operatorname{Re}\bigl(e^{i\theta}\langle a,\lambda\rangle\bigr)
 \ge\delta |a|.
\]
Thus, for each $c\in\C$, only finitely many monomials have $S$-weight $c$. Every eigenvector of $S$ in $\C[[z]]$ is therefore a polynomial.

On every finite quotient, the projection onto an eigenspace is a polynomial in $S$ and preserves the image of $\widehat I$. Passing to the limit shows that every weight component of an element of $\widehat I$ belongs to $\widehat I$. Let $J\subset\C[z]$ be generated by these homogeneous polynomial components. By noetherianity, $J$ has finitely many homogeneous generators $g_1,\ldots,g_r$. Finite sums of components approximate every element of $\widehat I$ in the $\m$-adic topology; since ideals in $\C[[z]]$ are closed,
\[
 \widehat I=(g_1,\ldots,g_r)\C[[z]].
\]

It remains to replace the complex weights by positive integers. The weight vectors $w$ for which every $g_j$ is homogeneous form the solution space $V\subset\C^N$ of the equations
\[
 \langle a-b,w\rangle=0,
\]
where $z^a,z^b$ occur in the same $g_j$. These equations have integer coefficients. The vector $e^{i\theta}\lambda$ lies in $V$, and since the equations are real so does its real part $\operatorname{Re}(e^{i\theta}\lambda)$, which has positive entries. Since $V\cap\R^N$ is defined over $\Q$, it contains a positive rational vector; clearing denominators gives positive integer weights.

The analytic germ defined by $g_1,\ldots,g_r$ is weighted homogeneous and formally isomorphic to $X$. The analytic isomorphism theorem following from Artin approximation \cite{Artin} makes these germs analytically isomorphic. Hence $X$ is weighted homogeneous.
\end{proof}

\begin{proof}[Proof of Theorem~\ref{thm:main}]
The zero-dimensional reduced germ is a point, for which the equivalences are immediate.

For a weighted homogeneous germ with weights $w_i>0$, the Euler field
\[
 E=\sum_iw_i z_i\partial_{z_i}
\]
is tangent to $X$, has positive spectrum, and points strictly outward along every small round link, since its radial component is a positive multiple of $\sum_iw_i|z_i|^2$. This gives \eqref{item:wh}$\Rightarrow$\eqref{item:single}. For the converse, apply Lemma~\ref{lem:spectrum} to the negative of the field in \eqref{item:single} to obtain \eqref{item:spectrum}. Lemma~\ref{lem:grading} then gives \eqref{item:wh}. On a connected link, the radial component of a transversal field has constant sign. Replacing the field by its negative if necessary gives strict outward pointing.
\end{proof}

\begin{corollary}\label{cor:all-links}
A reduced complex analytic germ is weighted homogeneous if and only if, in suitable local coordinates, it admits a holomorphic vector field pointing strictly outward along every sufficiently small link. If the small links are connected, it suffices that the field be transversal to every sufficiently small link.
\end{corollary}

\begin{proof}
The Euler field gives the forward implication. For the converse, let $v$ point strictly outward along every sufficiently small link. If $v(0)\ne0$, the real-time integral curve $\gamma$ through $0$, parametrized by $\dot\gamma=v(\gamma)$ in ambient coordinates, lies in $X$ and satisfies $\gamma(s)=s\,v(0)+O(s^2)$. Thus
\[
 \frac{d}{ds}\,\rho(\gamma(s))=2s\,|v(0)|^2+O(s^2)
\]
is negative for small $s<0$, contrary to strict outward pointing at $\gamma(s)$. Hence $v(0)=0$ and Theorem~\ref{thm:main} applies.

Suppose instead that the small links are connected and $v$ is transversal to each of them. In positive dimension, choose a nonconstant complex analytic arc $\alpha$ in $X$ through $0$. For small positive real $s$, the function $\rho(\alpha(s))$ is strictly increasing and takes all sufficiently small positive values. The radial component of $v$ has constant sign along the arc, and connectedness gives the same sign on every sufficiently small link. After changing the sign of $v$ if necessary, the field points strictly outward on all these links.
\end{proof}

\subsection*{AI declaration}
During the preparation of this manuscript, the authors used ChatGPT-6, OpenAI, and Claude Opus 5.5 (Anthropic), as aids in mathematical exploration, the discussion and refinement of arguments, and the improvement of exposition. All mathematical results and proofs have been critically examined and verified by the authors, who take full responsibility for the content of the paper.

\end{document}